\documentclass[11pt]{article}
\usepackage{times}
\usepackage{amsfonts}
\usepackage{amsmath}
\usepackage{amssymb}
\usepackage{latexsym,color}
\usepackage{amsbsy}
\usepackage{graphicx}

\newtheorem{thm}{Theorem}[section]
\newtheorem{con}{Conjecture}[section]
\newtheorem{lemma}[thm]{Lemma}
\newtheorem{cor}[thm]{Corollary}
\newtheorem{pro}[thm]{Proposition}
\newtheorem{example}[thm]{Example}
\newtheorem{definition}[thm]{Definition}

\newtheorem{remark}[thm]{Remark}
\newtheorem{Algorithm}[thm]{Algorithm}

\newcommand{\comment}[1]{}

\newcommand{\N}{{\mathbb N}}
\newcommand{\Bq}{B_q}
\newcommand{\Cq}{C_q}

\newcommand{\Fq}{{\mathbb F}_q}

\newcommand{\pr}{\parallel}

\newcommand{\kq}[1]{{{ {\bf{[#1]}} }}}
\newcommand{\qb}[2]{{{ {\bf{{#1}\brack {#2}}}}} }
\newcommand{\bin}[2]{{{ {#1}\choose {#2}}} }
\newcommand{\bfx}[1]{{{\bf{#1}}}}   
\newcommand{\oc}[1]{{{ \overline{#1} }}}

\newcommand{\ncom}{\newcommand}
\ncom{\ns}{\normalsize}
\ncom{\la}{\lambda}
\ncom{\bm}{\boldmath}
\ncom{\noi}{\noindent}
\ncom{\bq}{\begin{equation}}
\ncom{\eq}{\end{equation}}
\ncom{\beqn}{\begin{eqnarray*}}
\ncom{\eeqn}{\end{eqnarray*}}
\ncom{\ba}{\begin{array}}
\ncom{\ul}{\underline}

\ncom{\ea}{\end{array}}
\ncom{\beq}{\begin{eqnarray}}
\ncom{\eeq}{\end{eqnarray}}
\ncom{\nno}{\nonumber}
\ncom{\hs}{\mbox{\hspace{.25cm}}}
\ncom{\rar}{\rightarrow}
\ncom{\Rar}{\Rightarrow}
\ncom{\noin}{\noindent}
\ncom{\bc}{\begin{center}}
\ncom{\ec}{\end{center}}
\ncom{\sz}{\scriptsize}
\ncom{\fpd}{\Phi(\pi^{'})}
\ncom{\fp}{\Phi(\pi) }
\ncom{\nk}{\left< \begin{array}{c}
                       n\\k \end{array} \right>}
\ncom{\nd}{1^{'},2^{'},\cdots,n^{'}}
\ncom{\R}{I\!\!R}
\ncom{\de}{\bigtriangleup (F_{2n},\leq)}
\ncom{\del}{\bigtriangleup}
\ncom{\cov}{<\!\!\!\!\cdot }
\ncom{\bt}{\begin{thm}}
\ncom{\bcon}{\begin{con}}
\ncom{\et}{\end{thm}}
\ncom{\econ}{\end{con}}
\ncom{\bl}{\begin{lemma}}
\ncom{\el}{\end{lemma}}
\ncom{\bco}{\begin{cor}}
\ncom{\ds}{\displaystyle}
\ncom{\eco}{\end{cor}}
\ncom{\bp}{\begin{pro}}
\ncom{\ep}{\end{pro}}
\ncom{\bex}{\begin{example}}
\ncom{\eex}{\end{example}}
\ncom{\bd}{\begin{definition}}
\ncom{\ed}{\end{definition}}
\ncom{\brm}{\begin{remark}}
\ncom{\erm}{\end{remark}}
\ncom{\bal}{\begin{Algorithm}}
\ncom{\eal}{\end{Algorithm}}
\ncom{\ol}{\overline}

\ncom{\pf}{\noi {\bf Proof  }}
\ncom{\be}{\begin{enumerate}}
\ncom{\ee}{\end{enumerate}}
\ncom{\s}{\subset}
\ncom{\T}{{\cal T}}
\ncom{\B}{{\cal B}}
\ncom{\A}{{\cal A}}
\ncom{\Z}{{\cal Z}}

\title{\Large{{\bf The complexity of the $q$-analog of the $n$-cube }}}

\author{{ { Murali K. Srinivasan}} \\
{\em  \normalsize{Department of Mathematics}}\\
{\em  \normalsize{Indian Institute of Technology, Bombay}}\\
{\em  \normalsize{Powai, Mumbai 400076, INDIA}}\\
{\bf  \texttt{mks@math.iitb.ac.in}}\\
{\bf  \texttt{murali.k.srinivasan@gmail.com}}
\\{\small Mathematics Subject Classifications: 05C50, 05E25.}
}

\begin{document}
\date{}
\maketitle

\begin{abstract}

We present a positive, combinatorial, good formula for the complexity (=
number of spanning trees) of the $q$-analog of the $n$-cube. Our method also
yields the explicit block diagonalization of the commutant of the
$GL(n,\Fq)$ action on $\Bq(n)$, the set of all subspaces of an
$n$-dimensional vector space over $\Fq$. 

\end{abstract}

{{\bf  \section {  Introduction }}}  

This paper is a revised and expanded version of part of a 
previous paper by the author {\bf\cite{s2}}.
We begin with a discussion of our motivating problem. 
The number of spanning trees of a graph $G$ is called the {\em
complexity} of $G$ and denoted $c(G)$. The {\em hypercube} $C(n)$ 
is the graph whose vertex set is
the set $B(n)$ of all subsets of the $n$-set $[n]=\{1,2,\ldots ,n\}$ and 
where two
subsets $X,Y \in B(n)$ are connected by an edge iff 
$X\subseteq Y$ or $Y\subseteq X$, and $||X|-|Y|| = 1$.
A classical result states that the eigenvalues of
the Laplacian of $C(n)$ are $2k,\;k=0,\ldots ,n$, with respective
multiplicities ${n\choose k}$. It follows from the matrix tree theorem that
\beq \label{hc}
c(C(n)) &=& \frac{1}{2^n}
\left\{\prod_{k=1}^{n} (2k)^{n\choose k}\right\} \;=\; 
\prod_{k=2}^{n} (2k)^{n\choose k}. 
\eeq
We now define a $q$-analog of $C(n)$. 
Let $q$ be a prime power and  let $\Bq(n)$ denote 
the set of all subspaces of an $n$-dimensional vector space over the finite
field $\Fq$.
The {\em $q$-binomial coefficient}
$\qb{n}{k}$ denotes the number of $k$-dimensional subspaces in $\Bq(n)$.
The {\em $q$-analog $\Cq(n)$ of the hypercube} 
is the graph whose vertex set is $\Bq(n)$,
and where  
subspaces $X,Y\in \Bq(n)$ are connected by an edge iff 
$X\subseteq Y$ or $Y\subseteq X$, and $|\mbox{dim}(X) - \mbox{ dim}(Y)| = 1$.

The problem of finding a formula for $c(\Cq(n))$ is significantly more
involved than the classical (i.e., $q=1$) case. The underlying reason
seems to be that while the 
collection of all subsets forms an association scheme the collection
of all subspaces does not. 
To the best of our knowledge the eigenvalues of the Laplacian of $\Cq(n)$
are not known. In this paper we present 
a positive, combinatorial, good
formula for $c(\Cq(n))$. This stops well short of actually finding the
eigenvalues of the Laplacian but can be used to efficiently write down
$c(\Cq(n))$ for any given $n$.
Let us explain this. 
Consider an algorithm that, on input $n$, writes down the number
$c(C(n))$ as the output. Now this number is large, having
exponentially many bits in its binary representation and therefore, 
to write it
down explicitly will take time exponential in $n$. However, note that 
the binomial coefficients $\bin{n}{k}$ have atmost 
$n$ bits in their binary representation and
can be calculated in time polynomial in $n$ using the Pascal triangle. Thus 
formula
(\ref{hc}) above shows that, given $n$, we can write down the number
$c(C(n))$ in time polynomial in $n$ using product and 
exponential notation.
So we define a formula for $c(\Cq(n))$ to be {\em good} if, given $n$, we can
write down $c(\Cq(n))$ in time polynomial in $n$ 
using sum, product and exponential notation (treating $q$
symbolically).
The terms positive and combinatorial will be self
explanatory after the statement of Theorems \ref{stvs} and
\ref{qc} below.

Let us first  
reformulate the original formula (\ref{hc})  for $c(C(n))$ in
order to bring out the similarity with our formula for $c(\Cq(n))$.
Note that the following reformulation
is also a good formula for $c(C(n))$.

We have 
\beq \nonumber 
c(C(n)) &= & \frac{1}{2^n} \left\{\ds {\prod_{k=1}^n (2k)^{n\choose k}
}\right\}\\ 
\nonumber &=& \frac{1}{2^n} \left\{{\ds{\prod_{k=1}^{n} (2k) }}\right\}
\left\{{\ds{\prod_{k=1}^{\lfloor n/2 \rfloor}}}
\left({\ds{\prod_{j=k}^{n-k}}}
(2j)\right)^{{n\choose
k}-{n\choose{k-1}} } \right\}\\
   \label{hcmks}    &=&   n!\left\{
{\ds{\prod_{k=1}^{\lfloor n/2 \rfloor}}}
\left({\ds{\prod_{j=k}^{n-k}}}
(2j)\right)^{{n\choose
k}-{n\choose{k-1}} } \right\}
\eeq
To see the equivalence of the first and second lines above note that,
for $1\leq j\leq n/2$, the exponent of $2j$ in the numerator of the first line is
$n\choose j$ and in the numerator of the second line is also ${n\choose j}=1 + {n\choose 1}- {n\choose
0}+\cdots +{n\choose j}-{n\choose {j-1}}$. Since ${n\choose k} = {n\choose
{n-k}}$ the same conclusion holds for $n/2 \leq j \leq n$.

For $k\in \N = \{0,1,2,\ldots \}$ set $\kq{k}=1+q+q^2+\cdots +q^{k-1}$. Let $k,n \in \N$ with $k\leq n/2$. For
$k \leq j \leq n-k+1$, define polynomials $F_q(n,k,j)$ in $q$ 
using the following recursion: 
\beqn
F_q(n,k,n-k+1)&=&1 \\
F_q(n,k,n-k) &=& \kq{k} + \kq{n-k}
\eeqn
and, for $k\leq j < n-k$,
\beq \label{dr}
 F_q(n,k,j) &=&  
(\kq{j} + \kq{n-j})F_q(n,k,j+1) -
(q^k\kq{j+1-k}\kq{n-k-j})F_q(n,k,j+2).
\eeq
In Section 2 we prove the following formula for $c(\Cq(n))$, 
which is similar to 
formula (\ref{hcmks}) above, except
that the explicit term 
$\prod_{j=k}^{n-k}(2j)$ is replaced by a recursive calculation.
\bt \label{stvs} We have
\beqn
c(\Cq(n)) &=& F_q(n,0,1) \left\{
{\ds{\prod_{k=1}^{\lfloor n/2 \rfloor}}}
F_q(n,k,k)^{ \qb{n}{k}  - \qb{n}{k-1} } \right\},
\eeqn
where $F_q(n,0,1) = \kq{1}\kq{2}\cdots \kq{n}.$ 
\et
Given $n$, the polynomials $F_q(n,k,j)$ and $\qb{n}{k} - \qb{n}{k-1}$
can be efficiently calculated (in time polynomial in $n$) ,
using the recurrence (\ref{dr}) and the $q$-Pascal triangle respectively, 
and the formula above for $c(\Cq(n))$ is clearly 
good in the technical sense we
have defined. The following table, computed using Maple, 
gives the first five values of $c(\Cq(n))$.
\beqn
c(\Cq(1)) &= & 1\\
c(\Cq(2)) &= & \kq{2} 2^q \\
c(\Cq(3)) &= & \kq{2} \kq{3} (4+3q+q^2)^{q(1+q)}\\
c(\Cq(4)) &= & \kq{2} \kq{3} \kq{4}
          (8+12q+12q^2+10q^3+4q^4+2q^5)^{q(1+q+q^2)}\\
          &&\times (2+2q)^{q^2(q^2+1)}\\
c(\Cq(5)) &= &
\kq{2} \kq{3} \kq{4} \kq{5}
F_q(5,1,1)^{q(1+q)(1+q^2)}\\
&& \times F_q(5,2,2)^{q^2(1+q+q^2+q^3+q^4)}
\eeqn
where
$F_q(5,2,2) = 4+8q+7q^2+4q^3+q^4$ 
and 
$$F_q(5,1,1) =
16+36q+53q^2+65q^3+69q^4+58q^5+42q^6+26q^7+13q^8+5q^9+q^{10}.$$

The table above suggests that the formula for $c(\Cq(n))$ in Theorem
\ref{stvs} is positive, i.e., for $0\leq k\leq n/2$, both $\qb{n}{k}  -
\qb{n}{k-1}$ and $F_q(n,k,k)$ have nonnegative coefficients (as polynomials
in $q$).
A special case of a result of Butler {\bf\cite{b}} shows that indeed
$\qb{n}{k}  - \qb{n}{k-1},\;k\leq n/2$ has nonnegative
coefficients. In Theorem \ref{qc} of Section 3 we show that 
the polynomials $F_q(n,k,j)$ also have
nonnegative coefficients by giving an explicit positive combinatorial
formula for them. 

In {\bf\cite{s1}} we studied the singular values of the 
up operator on the Boolean
algebra and used this to explicitly block diagonalize the commutant of the
symmetric group action on $B(n)$. Theorem \ref{mt3} of
Section 2 gives the singular values of the up operator on the $q$-analog of the
Boolean algebra and we then 
use this result to
prove Theorem \ref{stvs}, which in essence is a block diagonalization of the
Laplacian of $\Cq(n)$. The Laplacian lies in 
the commutant of the $GL(n,\Fq)$ action on $\Bq(n)$
and
in Theorem \ref{mks} of Section 4 we generalize Theorem \ref{stvs} by   
explicitly block diagonalizing this commutant.

{\textcolor{black} {\bf \section {Singular values}}}

 A (finite) {\em graded poset} is a (finite) poset $P$ together with a
{\em rank function}
$r: P\rar \N$ such that if $p'$ covers $p$ in $P$ then
$r(p')=r(p)+1$. The {\em rank} of $P$ is $r(P)=\mbox{max}\{r(p): p\in P\}$
and,
for $i=0,1,\ldots ,r(P)$, $P_i$ denotes the set of elements of $P$ of rank
$i$. 

For a finite set $S$, let $V(S)$ denote the complex vector space with $S$ as
basis. Let $P$ be a graded poset with $n=r(P)$. Then we have
$ V(P)=V(P_0)\oplus V(P_1) \oplus \cdots \oplus V(P_n)$ (vector space direct
sum).
An element $v\in V(P)$ is {\em homogeneous} if $v\in V(P_i)$ for some $i$,
and if $v\not= 0$, we extend the notion of rank to nonzero homogeneous elements by writing
$r(v)=i$. 
The {\em
up operator}  $U:V(P)\rar V(P)$ is defined, for $p\in P$, by
$U(p)= \sum_{p'} p'$,
where the sum is over all $p'$ covering $p$.
A {\em symmetric Jordan chain} (SJC) in $V(P)$ is a sequence 
\beq \label{gjc}
&s=(v_k,\ldots ,v_{n-k}),\;\;\;k\leq n/2,&
\eeq 
of nonzero homogeneous elements of $V(P)$
such that $r(v_i)=i$ for $i=k,\ldots ,n-k$, $U(v_{i-1})=v_i$, for
$i=k+1,\ldots , n-k$, and $U(v_{n-k})=0$ (note that the
elements of this sequence are linearly independent, being nonzero and of
different ranks). We say that $s$ {\em
starts} at rank $k$ and {\em ends} at rank $n-k$. 
A {\em symmetric Jordan basis} (SJB) of $V(P)$ is a basis of $V(P)$
consisting of a disjoint union of symmetric Jordan chains 
in $V(P)$.

Let $\langle , \rangle$ denote the standard inner product on $V(P)$,
i.e.,
$\langle p,p' \rangle =
\delta (p,p')$ (Kronecker delta) for $p,p'\in P$.
The {\em length} $\sqrt{\langle v, v \rangle }$ of $v\in V(P)$ is denoted
$\pr v \pr$.

Suppose we have an orthogonal SJB $J(n)$ of $V(P)$.
Normalize the vectors in $J(n)$ to get an
orthonormal basis $J'(n)$.
Let $(v_k,\ldots ,v_{n-k})$ be a SJC in $J(n)$.
Put
$v_u' = \frac{v_u}{\pr v_u \pr}$ and $\alpha_u = \frac{\pr v_{u+1} \pr}{\pr
v_{u} \pr},\;k\leq u \leq n-k$ (we set $v_{k-1}'=v_{n-k+1}=0$). We have,
for $k\leq u \leq n-k$,
\beq \label{trick}
U(v_{u}')=\frac{U(v_{u})}{\pr v_{u} \pr}=\frac{v_{u+1}}{\pr v_{u}
\pr}=\alpha_{u} v_{u+1}'.&&
\eeq
Thus the matrix of $U$ wrt $J'(n)$ is in block diagonal form, with a block
corresponding to each (normalized) SJC in $J(n)$, and with
the block corresponding to $(v_k',\ldots ,v_{n-k}')$ above being a lower triangular
matrix with subdiagonal $(\alpha_k ,\ldots ,\alpha_{n-k-1})$ and $0$'s
elsewhere.

The {\em
down operator}  $D:V(P)\rar V(P)$ is defined, for $p\in P$, by
$D(p)= \sum_{p'} p'$,
where the sum is over all $p'$ covered by  $p$.
Note that the matrices, in the standard basis $P$,
of $U$ and $D$ are real and transposes of
each other. Since $J'(n)$ is orthonormal
with respect to the standard inner product, it follows that the matrices of
$U$ and $D$, in the basis $J'(n)$, must be adjoints of each other.
Thus, for $k-1\leq u \leq n-k-1$,  we must have (using (\ref{trick}) and the
previous paragraph),
\beq \label{trick1}
D(v_{u+1}')=\alpha_{u} v_{u}'.&&
\eeq
In particular, the subspace spanned by $\{v_k,\ldots
,v_{n-k}\}$ is closed under $U$ and $D$. We use this observation and
identities (\ref{trick}) and (\ref{trick1}) above without
explicit mention in a few places in Section 4.

Another useful observation is the following: take scalars $\beta_0 ,\beta_1
,\ldots ,\beta_n$ and define the operator $\gamma : V(P)\rar V(P)$ by
$\gamma (p) = \beta_{r(p)} p,\;p\in P$. Since each element of the 
SJC $(v_k,\ldots ,v_{n-k})$ is homogeneous, it follows from the
definition of $\gamma$ that the subspace spanned by $\{v_k ,\ldots
,v_{n-k}\}$
is closed under $U,D$ and $\gamma$.

The {\em Boolean algebra} is the graded poset of rank $n$  obtained by partially
ordering $B(n)$ by containment (with rank of a subset given by cardinality). 
The {\em $q$-analog of the Boolean algebra} is obtained by
partially ordering $\Bq(n)$ by inclusion. This gives a graded poset of rank
$n$ with rank of a subspace given by dimension. The following result is the
$q$-analog of a result about $B(n)$ proved in {\bf\cite{s1}}, which in turn
was motivated by Schrijver's fundamental paper {\bf\cite{s}}. 
 
\bt \label{mt3}
There exists a SJB $J(q,n)$ of $V(\Bq(n))$ such that

\noi (i) The elements of $J(q,n)$ are orthogonal with respect to
$\langle , \rangle$ (the standard inner product).

\noi (ii) {\em (Singular Values)} Let $0\leq k \leq n/2 $ and let
$(v_k,\ldots ,v_{n-k})$ be any SJC
in $J(q,n)$ starting at rank $k$ and ending at rank $n-k$. Then
we have, for $k\leq u < n-k$,
\beq \label{mkssv}
\frac{\pr v_{u+1} \pr}{\pr v_u \pr} & = &
\sqrt{q^k\kq{u+1-k}\kq{n-k-u}}
\eeq
\et
\pf We shall put together several 
standard results.
 
(i) The map $U^{n-2k}:V(\Bq(n)_k)\rar
V(\Bq(n)_{n-k}),\;0\leq k \leq n/2$ is well known to be bijective. 
It follows, using a standard Jordan canonical
form argument, that an SJB of $V(\Bq(n))$ exists.

(ii) Now we show existence of an orthogonal SJB. We use the action 
of the group $GL(n,\Fq)$ on $\Bq(n)$.
As is easily seen the
existence of an orthogonal SJB of $V(\Bq(n))$ (under the standard inner
product)
follows from facts (a)-(d) below  by an application of Schur's lemma:

\noi (a) Existence of some SJB of $V(\Bq(n))$.

\noi (b) $U$ is $GL(n,\Fq)$-linear.

\noi (c) For $0\leq k \leq n$, $V(\Bq(n)_k)$ is a mutiplicity free
$GL(n,\Fq)$-module (this is well known). 

\noi (d) For a finite group $G$, a $G$-invariant inner product on an
irreducible $G$-module is unique upto scalars.

(iii) Now we prove part (ii) of the Theorem. 
Define an operator $H: V(\Bq(n))\rar V(\Bq(n))$ by
$$H(X)= (\kq{k} - \kq{n-k})X,\;\; X\in \Bq(n)_k,\;0\leq k \leq n.$$ 
It is easy to check that $[U,D] = UD - DU = H$. To see this, fix $X\in
\Bq(n)_k$, and note that 
$UD(X) = \kq{k} X + \sum_Y Y$, where the sum is
over all $Y\in \Bq(n)_k$ with $\mbox{dim}(X\cap Y) = k-1$. 
Similarly, $DU(X) = \kq{n-k} X + \sum_Y Y$, where the sum is
over all $Y\in \Bq(n)_k$ with $\mbox{dim}(X\cap Y) = k-1$. 
Subtracting we get  $[U,D]=H$.

Let $J(q,n)$ be an orthogonal SJB of $V(\Bq(n))$ and let $(v_k,\ldots
,v_{n-k})$ be a SJC in $J(q,n)$ starting at rank $k$ and ending at rank $n-k$.
Put
$v_j' = \frac{v_j}{\pr v_j \pr}$ and $\alpha_j = \frac{\pr v_{j+1} \pr}{\pr
v_{j} \pr},\;k\leq j \leq n-k$. We have, 
from (\ref{trick}) and (\ref{trick1}),
$$U(v_j')=\alpha_j v_{j+1}',\;\;D(v_{j+1}')=\alpha_j v_j',\;\;k\leq j <
n-k.$$
We need to show that 
\beq \label{ind}
\alpha_j^2 &=& q^k \kq{j+1-k}\kq{n-k-j},\;\;k\leq j < n-k .
\eeq 
We show this by
induction on $j$. We have $DU = UD - H$. Now $DU(v_k')=\alpha_k D(v_{k+1}')
= \alpha_k^2 v_k'$ and $(UD-H)(v_k')= (\kq{n-k} - \kq{k}) v_k'$ (since
$D(v_k')=0$). Hence
$\alpha_k^2 = \kq{n-k} - \kq{k} = q^k \kq{n-2k}$. Thus (\ref{ind}) holds for
$j=k$.

As in the previous paragraph $DU(v_j')=\alpha_j^2 v_j'$ and
$(UD-H)(v_j')=(\alpha_{j-1}^2 + \kq{n-j} - \kq{j}) v_j'$.
By induction, we may assume $\alpha_{j-1}^2 = q^k \kq{j-k}\kq{n-k-j+1}$.
Thus we see that $\alpha_j^2$ is    
\beqn
&=& q^k \kq{j-k}\kq{n-k-j+1} + \kq{n-j} - \kq{j}\\
&=& q^k \left\{ (\kq{j+1-k} - q^{j-k})(\kq{n-k-j} + q^{n-k-j})\right\}
+ \kq{n-j} - \kq{j}\\
&=& q^k \left\{ \kq{j+1-k}\kq{n-k-j} + q^{n-k-j}\kq{j+1-k}
 -q^{j-k}\kq{n-k-j} - q^{n-2k}\right\} \\
&& + \kq{n-j} - \kq{j}\\
&=& q^k \kq{j+1-k}\kq{n-k-j} + q^{n-j}\kq{j+1-k} - q^j\kq{n-k-j} - q^{n-k}\\
&&     + \kq{n-j} - \kq{j}\\
&=& q^k \kq{j+1-k}\kq{n-k-j} + \kq{n+1-k} - \kq{n-j} - \kq{n+1-k} + \kq{j}\\
  &&   + \kq{n-j} - \kq{j}\\
&=& q^k \kq{j+1-k}\kq{n-k-j},
\eeqn
completing the proof.$\Box$

{\bf Remark} The proof in {\bf\cite{s1}} of the $q=1$ case of Theorem
\ref{mt3} was constructive, giving a simple algorithm to explicitly write
down an orthogonal SJB of $V(B(n))$. It would be interesting to give an 
explicit construction of an orthogonal SJB of $V(\Bq(n))$.

For $0\leq k \leq n/2$, define a real, symmetric, tridiagonal matrix
$N=N(k,n-k,n)$ of size $n-2k+1$,
with rows and columns indexed by the set $\{k,k+1,\ldots ,n-k\}$, and with
entries given as follows: 
for $k\leq i,j \leq n-k$ define
\beqn
N(i,j) &=& \left\{ \ba{ll}
      -\sqrt{q^k\kq{j-k}\kq{n-k-j+1}} & \mbox{if $i=j-1$}\\
                                 & \\ 
      \kq{j} + \kq{n-j}  & \mbox{if $i=j$}\\
                    & \\
      -\sqrt{q^k\kq{j+1-k}\kq{n-k-j}} & \mbox{if $i=j+1$}\\
        & \\
      0 & \mbox{if $|i-j|\geq 2$}
                   \ea \right.
\eeqn
For $0\leq k \leq n/2 $ and $k\leq j \leq n-k+1$ define
$N_j=N_j(k,n-k,n)$
to be the principal submatrix of $N=N(k,n-k,n)$ indexed by the rows and
columns in the set $\{j,j+1,\ldots ,n-k\}$. Thus, $N_k = N$ and $N_{n-k+1}$
is the empty matrix, which by convention has determinant 1.

\bl \label{det}
 
For $0\leq k \leq n/2 $ and $k\leq j \leq n-k+1$ we have 
 
(i) $F_q(n,k,j) = \mbox{\em det}(N_j(k,n-k,n)).$

(ii) $F_q(n,0,j) = \kq{n}\kq{n-1}\cdots \kq{j}$. 

In particular, $\mbox{\em det}(N_1(0,n,n)) = \kq{n}\kq{n-1}\cdots \kq{1}$. 
\el
\pf (i) By (reverse) induction on $j$. The base cases $j=n-k+1,n-k$ are 
clear and the general case follows by expanding the determinant of $N_j$
along its first column.

(ii) By (reverse) induction on $j$. The base cases $j=n+1,n$ are clear. By
induction and the defining recurrence for $F_q(n,k,j)$ we have 
\beqn
F_q(n,0,j) & = & (\kq{j} + \kq{n-j})F_q(n,0,j+1) -
(\kq{j+1}\kq{n-j})F_q(n,0,j+2)\\
& = & (\kq{j} + \kq{n-j})\kq{n}\cdots \kq{j+1} -
(\kq{j+1}\kq{n-j})\kq{n}\cdots \kq{j+2}\\
&=& \kq{n}\cdots \kq{j}, 
\eeqn
completing the proof.$\Box$

We now prove our first formula for $c(\Cq(n))$.

{\bf Proof} {\em (of Theorem \ref{stvs})} 

The degree of a vertex $X$
of $\Cq(n)$ is $\kq{k} + \kq{n-k}$, where $k=\mbox{dim}(X)$. 
Define an operator $deg: V(\Bq(n))\rar
V(\Bq(n))$ by
$$ deg(X) = (\kq{k}+\kq{n-k})X,\;X\in \Bq(n)_k.
$$
We can now write the Laplacian $L: V(\Bq(n))\rar V(\Bq(n))$ of $\Cq(n)$  as
$ L = deg - U - D,
$
where $U,D$ are the up and down operators on $V(\Bq(n))$.

Let $J(q,n)$ be a SJB of $V(\Bq(n))$ satisfying the conditions
of Theorem \ref{mt3}. Normalize $J(q,n)$ to get an orthonormal basis
$J'(q,n)$. Since the vertex degrees are constant on $\Bq(n)_k$ it follows 
that the subspace spanned by each SJC in $J(q,n)$ is closed under $L$. Using
part (ii) of Theorem \ref{mt3} we can write down the matrix of
$L$ in the basis $J'(q,n)$.

Let $0\leq k \leq n/2$. Let $(w_k,\ldots
,w_{n-k})$ be a SJC in $J(q,n)$ starting at rank $k$.
Set $v_i = \frac{w_{i}}{\pr w_{i} \pr},\;k\leq i \leq n-k$. Let $W$ be the
subspace spanned by $\{v_k , \ldots ,v_{n-k}\}$. Then $W$ is invariant
under $L$.

It follows from Theorem \ref{mt3} that
$N(k,n-k,n)$ is 
the matrix of $L:W\rar W$ with
respect to the (ordered) basis $\{v_k ,\ldots ,v_{n-k}\}$ 
(we take coordinate vectors with respect to a basis as column vectors).
Thus the matrix of
$L$ with respect to (a suitable ordering of) $J'(q,n)$ is in block diagonal form,
with
blocks $N(k,n-k,n)$, for all $0\leq k \leq n/2$, and 
each such block is repeated 
$\qb{n}{k} - \qb{n}{k-1}$ times. The number of distinct blocks is $1 +
\lfloor n/2 \rfloor$.

The unique element in $J'(q,n)$ of rank 0 is the vector ${\bf 0}$ (here ${\bf
0}$ is the zero subspace).

Let ${\cal M}$ denote the matrix of the Laplacian of $\Cq(n)$ in the standard
basis $\Bq(n)$ and let ${\cal M}'$ be obtained from ${\cal M}$ by removing the row and column
corresponding to vertex ${\bf 0}$.
From the matrix tree
theorem we have  
$c(\Cq(n))=\mbox{ det}({\cal M}')$.
A little reflection shows that, by changing bases from 
$\Bq(n) - \{{\bf 0}\}$ to $J'(q,n) - \{{\bf 0}\}$, 
${\cal M}'$ block diagonalizes
with a block $N_1(0,n,n)$ of multiplicity 1 and blocks $N(k,n-k,n)$,  
$1\leq k \leq n/2$, of multiplicity $\qb{n}{k} - \qb{n}{k-1}$. The result
now follows from Lemma \ref{det}.$\Box$ 

{\bf Remark}
Since the subdiagonal entries of $N(k,n-k,n),\;0\leq k
\leq n/2$ are nonzero it easily follows that any eigenspace will have
dimension 1 and thus $N(k,n-k,n)$ has $n-2k+1$ distinct eigenvalues. 
Data suggest that 
$N(k,n-k,n)$ and $N(l,n-l,n)$, $k\not=l$ do not have any eigenvalue in
common. In other words, the Laplacian of $\Cq(n)$ seems to have
$\sum_{k=0}^{\lfloor n/2 \rfloor} (n-2k+1) = (\lfloor n/2 \rfloor + 1)(\lceil n/2 \rceil +
1)$ distinct eigenvalues, with each of the $n-2k+1$ 
eigenvalues of $N(k,n-k,n)$ having multiplicity
$\qb{n}{k} - \qb{n}{k-1},\;k=0,1,\ldots , \lfloor n/2 \rfloor$. Is this
true and can it be proved without explicitly writing down the eigenvalues.

{\textcolor{black} {\bf \section {Positivity}}}

In this section we define certain combinatorial objects and a generating function
based on them with the property that an
appropriate positive specialization satisfies the recurrence (\ref{dr}).

Let $[\oc{n}]=\{\oc{1},\oc{2},\ldots ,\oc{n}\}$ and 
consider the set $[n,\oc{n}]= [n] \cup [\oc{n}]$ of $2n$ elements. 
We are going to recursively define a set $S(n)$
of certain subsets of $[n,\oc{n}]$. The cardinality of an element of $S(n)$ will be
between $0$ and $n$ (inclusive) and will have the same parity as $n$. Define
$S(0) = \{ \emptyset \}$, $S(1)=\{ \{1\}, \{\oc{1}\} \}$ and, for $n\geq 1$,
\beqn \label{co}
S(n+1) &=&  
\{ X\cup\{n+1\} : X\in S(n)\} \cup \{ X\cup \{\oc{n+1}\} : X\in
S(n),\; n \not\in X\} \cup S(n-1).
\eeqn
It is easy to show by induction that $|S(n)| = 2^n$ and that the number of
elements of $S(n)$ not containing $n$ is $2^{n-1}$. Let $\bfx{x} =
(x_1,x_2,\ldots ),\; \bfx{y} = (y_1,y_2,\ldots )$, and $z$ be
indeterminates. For $n\geq 0$ define the following polynomial
$$ P(n,\bfx{x},\bfx{y},z) = \sum_{X\in S(n)} \left( \prod_{i\in X\cap [n]}
x_i \right) \left( \prod_{\oc{i}\in X \cap [\oc{n}]} y_i \right)
z^{\frac{n-|X|}{2}}.$$
{\bf Example} We have 
\beqn S(2) &=& \{ \emptyset , \{1,2\}, \{\oc{1},2\} ,
\{\oc{1},\oc{2}\} \},\\
S(3)&=&\{ \{1\}, \{\oc{1}\}, \{3\}, \{\oc{3}\}, \{1,2,3\}, \{\oc{1},2,3\},
\{\oc{1},\oc{2},3\}, \{\oc{1},\oc{2},\oc{3}\} \}.
\eeqn
Thus $P(2,\bfx{x},\bfx{y},z)= z+ (x_1x_2 + y_1x_2 + y_1y_2)$ and 
$$P(3,\bfx{x},\bfx{y},z) = (x_1 + y_1 + x_3 + y_3)z + (x_1x_2x_3 + y_1x_2x_3
+ y_1y_2x_3 + y_1y_2y_3).$$
The recursive structure of $S(n)$ yields the following recurrence for the
polynomials $P$.
\bt \label{xy}
We have
$$P(n+1,\bfx{x},\bfx{y},z) = (x_{n+1} + y_{n+1}) P(n,\bfx{x},\bfx{y},z)
                             - (x_n y_{n+1}
                             - z)P(n-1,\bfx{x},\bfx{y},z),\;\;n\geq1.$$
\et
\pf Let $X\in S(n+1)$. In the
expansion of the lhs $P(n+1,\bfx{x},\bfx{y},z)$,
consider the term corresponding to $X$: 
$$ \left( \prod_{i\in X\cap [n]}
x_i \right) \left( \prod_{\oc{i}\in X \cap [\oc{n}]} y_i \right)
z^{\frac{n+1-|X|}{2}}.$$
We consider three cases:

\noi (i) $n+1\in X$: the term above will appear exactly once in
$x_{n+1}P(n,\bfx{x},\bfx{y},z)$.

\noi (ii) $\oc{n+1}\in X$: the term above will appear exactly once 
in $y_{n+1} P(n,\bfx{x},\bfx{y},z) - x_ny_{n+1} P(n-1,\bfx{x},\bfx{y},z)$.

\noi (iii) $X\in S(n-1)$: The term above will appear exactly once in
$zP(n-1,\bfx{x},\bfx{y},z)$.

It is clear that there are no other terms corresponding to $X$ 
on the rhs. The result follows.
$\Box$

Given $n,k\in \N$ with $k\leq n/2$ define
$$
d_q(n,k) = (\kq{n-k},\kq{n-k-1},\ldots ,\kq{k},0,0,\ldots ),\;\;
e_q(n,k) = (\kq{k},\kq{k+1},\ldots ,\kq{n-k},0,0,\ldots ).
$$
We now prove the nonnegativity of the coefficients of $F_q(n,k,j)$.
\bt \label{qc}
Let $n,k\in \N$ with $k\leq n/2$ and let $k\leq j \leq n-k+1$. Then
$$F_q(n,k,j) =
P(n-k-j+1,\;d_q(n,k),\;e_q(n,k),\;\kq{k}\kq{n-k+1}).$$
\et
\pf The result is clearly true for $j=n-k+1$ and
$j=n-k$. Now note the following alternate
expression for the singular values:
\beqn
q^k \kq{j+1-k}\kq{n-k-j}&=&(\kq{j+1}-\kq{k})(\kq{n-j}-q^{n-k-j}\kq{k})\\
&=&\kq{j+1}\kq{n-j} - \kq{k}(\kq{n-j} + q^{n-k-j}\kq{j+1} -
q^{n-k-j}\kq{k})\\
&=&\kq{j+1}\kq{n-j} - \kq{k}\kq{n-k+1}
\eeqn 
It now follows from Theorem \ref{xy} that
$P(n-k-j+1,\;d_q(n,k),\;e_q(n,k),\;\kq{k}\kq{n-k+1})$ satisfies the
same recurrence as $F_q(n,k,j)$. The result follows. $\Box$

{\bf Example} 
Let $n=3$ and $k=0$. Then 
$$d_q(3,0) = (\kq{3},\kq{2},\kq{1},\kq{0},0,0,\ldots ),\;\;
e_q(3,0) =  (\kq{0},\kq{1},\kq{2},\kq{3},0,0,\ldots ).$$
Thus
$F_q(3,0,1)=P(3,d_q(3,0),e_q(3,0),0) = \kq{3}\kq{2}\kq{1}$.

Now let $n=3$ and $k=1$. Then 
$$d_q(3,1) = (\kq{2},\kq{1},0,0,\ldots ),\;\;
e_q(3,1) =  (\kq{1},\kq{2},0,0,\ldots ).$$
Thus
$F_q(3,1,1)=P(2,d_q(3,1),e_q(3,1),\kq{1}\kq{3}) = \kq{1}\kq{3} +
\kq{2}\kq{1} + \kq{1}\kq{1} + \kq{1}\kq{2} = 4+3q+q^2$, agreeing with the
formula given in the introduction.

Taking $d(n,k)=(n-k,n-k-1,\ldots ,k,0,0,\ldots)$,
$\;e(n,k)=(k,k+1,\ldots ,n-k,0,0,\ldots)$,
substituting $q=1$ in the formula above and comparing with (\ref{hcmks}) we
get the following identity 
\beqn 
P(n-2k+1,\;d(n,k),\;e(n,k),\;k(n-k+1))&=& 
2^{n-2k+1}\cdot k\cdot (k+1) \cdots (n-k).\eeqn

{\textcolor{black} {\bf \section {Explicit block diagonalization}}}

The group $G=GL(n,\Fq)$ 
has a rank and order preserving action on the graded poset $\Bq(n)$. 
Note that the Laplacian of $\Cq(n)$ belongs to the algebra $\mbox{End}_G
(V(\Bq(n)))$. In this section we generalize Theorem \ref{stvs} by 
explicitly block diagonalizing $\mbox{End}_G (V(\Bq(n)))$. 
Our method is the $q$-analog of the
method used in {\bf\cite{s1}} to explicitly block diagonalize the
commutant of the symmetric group action on $B(n)$.

We represent
elements of $\mbox{End}(V(\Bq(n)))$ (in the standard basis) as $\Bq(n)\times
\Bq(n)$ matrices (we think of elements of $V(\Bq(n))$ as column vectors with
coordinates indexed by $\Bq(n)$). For
$X,Y\in \Bq(n)$, the entry in row $X$, column $Y$
of a matrix $M$ will be denoted $M(X,Y)$. 
The matrix corresponding to $f\in \mbox{End}(V(\Bq(n)))$ is denoted
$M_f$. 
Set  
$\A(q,n) = \{ M_f : f \in \mbox{End}_G (V(\Bq(n)))\}.$
Then $\A(q,n)$ is a $*$-algebra of matrices.

Let $f\in \mbox{End}(V(\Bq(n)))$ and $g\in G$. Then
$$f(g(Y)) = \sum_{X}M_f(X,g(Y))X \mbox{ and }
g(f(Y)) = \sum_{X}M_f(X,Y)g(X).$$   
It follows that $f$
is $G$-linear if and only if
$M_f(X,Y) = M_f(g(X),g(Y))$,
for all $X,Y\in \Bq(n),\;g\in G$,
i.e., $M_f$ is constant on the orbits of the $G$-action on $\Bq(n) \times
\Bq(n)$. 

For $0\leq i,j,t \leq n$ let $M^t_{i,j}$ be the 
$\Bq(n)\times \Bq(n)$ matrix given by   
$$M^t_{i,j}(X,Y) = \left\{ \ba{ll}
                            1 & \mbox{if } \mbox{dim}(X)=i,\;
                            \mbox{dim}(Y)=j,\;\mbox{dim}(X\cap Y)=t \\
                            0 & \mbox{otherwise}
                            \ea
                    \right.
$$
Now $(X,Y),(X',Y')\in 
\Bq(n) \times \Bq(n)$ are in the same $G$-orbit iff
$\mbox{dim}(X)=\mbox{ dim}(X'), \mbox{ dim}(Y)=\mbox{ dim}(Y')$, and 
$\mbox{dim}(X\cap Y)=\mbox{ dim}(X'\cap Y')$.
It follows that 
$$\{M^t_{i,j}\;|\; i-t+t+j-t \leq n,\;i-t,t,j-t \geq 0 \}$$
is a basis of $\A(q,n)$ and its cardinality is $\bin{n+3}{3}$.

Fix $i,j\in\{0,\ldots, n\}$. 
Then we have
$$M_{i,t}^tM_{t,j}^t = \sum_{u=0}^n \qb{u}{t}M_{i,j}^u,\;\;\;\;\;
t=0,\ldots ,n,$$
since the entry of the lhs in row $X$, col $Y$ with $\mbox{dim}(X)=i, \mbox{
dim}(Y)=j$ is equal
to the number of common subspaces of $X$ and $Y$ of size $t$. 
Apply $q$-binomial inversion (see Exercise 2.47 in {\bf\cite{a}}) to get
\beq \label{qbi}
M_{i,j}^t= \sum_{u=0}^n (-1)^{u-t}q^{\bin{u-t}{2}}\qb{u}{t}
M_{i,u}^uM_{u,j}^u,\;\;\;\;\;t=0,\ldots ,n.
\eeq

For the rest of this section
set $m=\lfloor n/2 \rfloor$, and $p_k 
=n-2k+1,\; q_k = {\qb{n}{k}}-{\qb{n}{k-1}},\; k=0,\ldots ,m$.
Note that
\beq \label{di}
\sum_{k=0}^m p_k^2 = {{n+3}\choose 3},
\eeq
since both sides are polynomials in $l$ (treating the cases $n=2l$ and
$n=2l+1$ separately) of degree 3 and agree for $l=0,1,2,3$.

Consider an orthogonal SJB $J(q,n)$ of $V(\Bq(n))$ 
satisfying the conditions of Theorem \ref{mt3}. Normalize $J(q,n)$ to get an
orthonormal basis $J'(q,n)$. Let $N(n)$ be the square
$\Bq(n) \times J'(q,n)$ matrix, where, for $v\in J'(q,n)$, the column of
$N(n)$ indexed by $v$ is the coordinate vector of $v$ (in the standard
coordinates $\Bq(n)$). 
By Theorem \ref{mt3}(i), $N(n)$ is unitary. 
Since the action of $M_{i,u}^u$ on
$V(\Bq(n)_u)$ is $\frac{1}{\kq{i-u}\kq{i-u-1}\cdots \kq{1}}$ 
times the action of $U^{i-u}$ on
$V(\Bq(n)_u)$, it follows  
by Theorem \ref{mt3}(ii) and identities (\ref{qbi}), (\ref{di}) above
that $N(n)^* \A(q,n) N(n)$ consists of all $J'(q,n) \times J'(q,n)$
block diagonal matrices with a block
corresponding to each (normalized) SJC in $J(q,n)$ and any two SJC's
starting and ending at the same rank giving rise to identical blocks.
So there are $q_k$ identical
blocks of size $p_k$, for $k=0,\ldots ,m$. It will be convenient to reindex
the rows and columns of a block corresponding to a SJC starting at rank $k$
and ending at rank $n-k$ by the set $\{k,k+1,\ldots ,n-k\}$.

Define a map (below $\mbox{Mat}(n\times n)$ denotes
the algebra of complex $n\times n$ matrices)
$$\Phi: \A(q,n) \cong \bigoplus_{k=0}^m \mbox{Mat}(p_k\times p_k),$$
by conjugating with
$N(n)$ followed by dropping duplicate blocks. We now write down
the image $\Phi(M_{i,j}^t)$.

For $i,j,k,t \in \{0,\ldots ,n\}$ define
\beqn
\beta_{i,j,k}^{n,t}(q) & = &\sum_{u=0}^n (-1)^{u-t} \;q^{\bin{u-t}{2}-ku}\;
\qb{u}{t}\qb{n-2k}{u-k}\qb{n-k-u}{i-u}\qb{n-k-u}{j-u}.
\eeqn
For $0\leq k \leq m$ and $k\leq i,j \leq n-k$, define $E_{i,j,k}$ to be the
$p_k \times p_k$ matrix, with rows and columns indexed by $\{k,k+1,\ldots  
,n-k\}$, and with entry in row $i$ and column $j$ equal to 1 and all other 
entries 0.

In the proof of the following result we will need another alternate expression for the
singular values:
\beq \label{asv}
\sqrt{q^k \kq{u+1-k}\kq{n-k-u}} &=&
q^{\frac{k}{2}}\;\kq{n-k-u}{\qb{n-2k}{u-k}}^{\frac{1}{2}}
{\qb{n-2k}{u+1-k}}^{- \frac{1}{2}}.
\eeq
\bt  \label{mks}
Let $i,j,t\in \{0,\ldots ,n\}$. Write
$$\Phi (M_{i,j}^t) = (N_0,\ldots ,N_m),$$
where, for
$k=0,\ldots , m$, the rows and columns of $N_k$
are indexed by $\{k,k+1,\ldots ,n-k\}$. 
Then, for $0\leq k \leq m$,
$$N_k = \left\{ \ba{ll}
                q^{\frac{k(i+j)}{2}}\;
                {\qb{n-2k}{i-k}}^{-\frac{1}{2}}
                {\qb{n-2k}{j-k}}^{-\frac{1}{2}}\;
                \beta_{i,j,k}^{n,t}(q)\;
                E_{i,j,k} & \mbox{if } k\leq i,j \leq n-k \\
                  0    & \mbox{otherwise}
                 \ea
         \right.
$$ 
\et

\pf Fix $0\leq k \leq m$. If both $i,j$ are not elements of $\{k,\ldots ,
n-k\}$ then clearly $N_k = 0$. So we may assume $k\leq i,j \leq n-k$.
Clearly, $N_k = \lambda E_{i,j,k}$ for some $\lambda$. We now find   
$\lambda = N_k(i,j)$.

Let $u\in \{0,\ldots ,n\}$.
Write 
$\Phi (M_{i,u}^u) = (A_0^u,\ldots ,A_m^u)$.
We claim that 
$$A_k^u = \left\{ \ba{ll}
q^{\frac{k(i-u)}{2}}\;{\qb{n-k-u}{i-u}}{\qb{n-2k}{u-k}}^{\frac{1}{2}}
                {\qb{n-2k}{i-k}}^{-\frac{1}{2}}\;
                E_{i,u,k} & \mbox{if } k\leq u \leq n-k \\
                  0    & \mbox{otherwise}
                 \ea
         \right.
$$ 
The otherwise part of the claim is clear. If $k\leq u \leq n-k$ and $i < u$
then we have $A_k^u = 0$. This also follows from the rhs since the $q$-binomial
coefficient $\qb{a}{b}$ is 0 for $b< 0$. So we may assume that $k\leq u   
\leq n-k$ and $i\geq u$. Clearly, in this case we have $A_k^u = \alpha     
E_{i,u,k}$, for some $\alpha$. We now determine $\alpha= A^u_k(i,u)$. We   
have using Theorem \ref{mt3}(ii) and the expression (\ref{asv}) 
\beqn
A_k^u(i,u) &=&  \frac{\prod_{w=u}^{i-1}\left\{ q^{\frac{k}{2}}\;
\kq{n-k-w}{\qb{n-2k}{w-k}}^{\frac{1}{2}}   
{\qb{n-2k}{w+1-k}}^{-\frac{1}{2}}\right\}}{\kq{i-u}\kq{i-u-1}\cdots
\kq{1}}\\&&\\
&=&
q^{\frac{k(i-u)}{2}}\;{\qb{n-k-u}{i-u}}{\qb{n-2k}{u-k}}^{\frac{1}{2}}
                {\qb{n-2k}{i-k}}^{-\frac{1}{2}}.
\eeqn
  
Similarly, if we write 
$\Phi (M_{u,j}^u) = (B_0^u,\ldots ,B_m^u)$,
then we have
$$B_k^u = \left\{ \ba{ll}
q^{\frac{k(j-u)}{2}}\;{\qb{n-k-u}{j-u}}{\qb{n-2k}{u-k}}^{\frac{1}{2}}
                {\qb{n-2k}{j-k}}^{-\frac{1}{2}}\;
                E_{u,j,k} & \mbox{if } k\leq u \leq n-k \\
                  0    & \mbox{otherwise}
                 \ea
         \right.
$$
So from (\ref{qbi}) we have 
$
N_k  =  \sum_{u=0}^n (-1)^{u-t} \;q^{\bin{u-t}{2}}\;{\qb{u}{t}} A_k^u B_k^u  
 =  \sum_{u=k}^{n-k} (-1)^{u-t} \;q^{\bin{u-t}{2}}\;{\qb{u}{t}} A_k^u B_k^u. 
$
Thus
\beqn
\lefteqn{N_k(i,j)}\\
 &=& \sum_{u=k}^{n-k} (-1)^{u-t} \;q^{\bin{u-t}{2}}\;{\qb{u}{t}}
\left\{ \sum_{l=k}^{n-k} A_k^u(i,l) B_k^u(l,j)\right\}\\ && \\
 &=& \sum_{u=k}^{n-k} (-1)^{u-t} \;q^{\bin{u-t}{2}}\;{\qb{u}{t}}
     A_k^u(i,u) B_k^u(u,j)\\ &&\\
 &=& \sum_{u=k}^{n-k} (-1)^{u-t} \;q^{\bin{u-t}{2}}\;{\qb{u}{t}}
\;q^{\frac{k(i-u)}{2}}\;{\qb{n-k-u}{i-u}}{\qb{n-2k}{u-k}}^{\frac{1}{2}}
                {\qb{n-2k}{i-k}}^{-\frac{1}{2}}\\ &&\\
&& \times q^{\frac{k(j-u)}{2}}\;{\qb{n-k-u}{j-u}}{\qb{n-2k}{u-k}}^{\frac{1}{2}}
                {\qb{n-2k}{j-k}}^{-\frac{1}{2}}\\  &&\\
& = & q^{\frac{k(i+j)}{2}}\;{\qb{n-2k}{i-k}}^{-\frac{1}{2}}
      {\qb{n-2k}{j-k}}^{-\frac{1}{2}}
      \left\{\sum_{u=0}^{n} (-1)^{u-t} \;q^{\bin{u-t}{2}-ku}\;{\qb{u}{t}}
{\qb{n-k-u}{i-u}}\right.\\ &&\\
&& \times \left.
{\qb{n-k-u}{j-u}}{\qb{n-2k}{u-k}}\right\},
\eeqn
completing the proof.$\Box$

{\bf Remark} Let $0\leq i \leq n/2$, $0\leq k,t \leq i$ and write
$\Phi (M_{i,i}^t) = (N'_0,\ldots ,N'_m).$
Substituting $j=i$ in Theorem \ref{mks} and noting that 
${\qb{n-2k}{i-k}}^{-1}{\qb{n-2k}{u-k}}{\qb{n-k-u}{i-u}} = {\qb{i-k}{i-u}}$
we see that $N'_k = \tau_{i,t,k} E_{i,i,k}$ where
$$\tau_{i,t,k} = \sum_{u=0}^n
(-1)^{u-t}\;q^{\bin{u-t}{2}+k(i-u)}\;\qb{u}{t}\qb{n-k-u}{i-u}\qb{i-k}{i-u}.$$
The $\tau_{i,t,k}$ are the eigenvalues of the $q$-Johnson
scheme of $i$-dimensional subspaces (see {\bf\cite{bi}}).

\begin{center} {\bf \Large{Acknowledgement}}
\end{center}
I thank Professors Alexander Schrijver, Navin Singhi, and Thomas Zaslavsky
for their comments on an earlier version of this paper. 
I am grateful to Sivaramakrishnan Sivasubramanian for several useful
discussions about the positivity of the polynomials $F_q(n,k,j)$.

\end{document}